\documentclass[a4paper,12pt,reqno]{amsart}
\usepackage[utf8]{inputenc}
\usepackage{csquotes}
\usepackage{amssymb}
\usepackage{amsmath}
\usepackage{makecell}
\usepackage{multirow}
\usepackage{hhline}
\usepackage{amsthm}
\usepackage{amsfonts}
\usepackage{mathtools}
\usepackage{todonotes}
\usepackage{comment}
\usepackage{mathrsfs}
\usepackage{enumitem}
\usepackage[all,cmtip]{xy}
\usepackage{tikz-cd}
\usepackage[hyphens]{url}
\usepackage{hyperref}
\usepackage{cleveref}
\usepackage[backend=biber, style=alphabetic]{biblatex}
\usepackage[british]{babel}
\usepackage[margin=6em]{geometry}
\usepackage{eqparbox}

\usepackage{pbox}

\newcommand{\set}[1]{\left\lbrace #1 \right\rbrace}

\newcommand{\field}[1]{\mathbb{#1}}  
\newcommand{\Q}{\field{Q}} 
\newcommand{\R}{\field{R}} 
\newcommand{\Z}{\field{Z}} 
\newcommand{\F}{\field{F}} 

\renewcommand{\P}{\field{P}}
\newcommand{\PP}{\field{P}}

\newcommand{\OO}{\mathcal{O}}

\DeclareMathOperator{\Div}{Div}
\DeclareMathOperator{\ddiv}{div}

\DeclareMathOperator{\Cl}{Cl}

\DeclareMathOperator{\Pic}{Pic}
\DeclareMathOperator{\Spec}{Spec}

\DeclareMathOperator{\Gal}{Gal}

\DeclareMathOperator{\Support}{Support}

\newtheorem{lemma}{Lemma}
\newtheorem{theorem}[lemma]{Theorem}
\newtheorem{proposition}[lemma]{Proposition}
\newtheorem{corollary}[lemma]{Corollary}

\theoremstyle{definition}
\newtheorem{definition}[lemma]{Definition}
\newtheorem{example}[lemma]{Example}

\newtheorem{question}[lemma]{Question}
\newtheorem{remark}[lemma]{Remark}

\numberwithin{lemma}{section}
\numberwithin{equation}{section} 
\numberwithin{figure}{section}

\newcommand{\githubbare}[1]{\href{https://github.com/koffie/mdmagma/blob/v0.2.1/computations/derickx-X1_16_classgroups/#1}{\path{#1}}}
\newcommand{\github}[2]{[\githubbare{#1},#2]}

\title{Class groups of imaginary quadratic points on $X_1(16)$}
\author{Maarten Derickx}
\date{July 2024}

\begin{document}
	
	\begin{abstract}
		The main goal is to show that if $K \ncong \Q(\sqrt{-15})$ is an imaginary quadratic field and $E$ is an elliptic curve over $K$ with a torsion point of order 16, then the class number of $K$ is divisible by 10. This gives an affirmative answer to a question by Krumm. This is done by setting up a more general framework for studying divisibility of class groups of imaginary quadratic points on hyper-elliptic curves.
	\end{abstract}

	\maketitle
	
	
	\section{Introduction}
	David Krumm studied quadratic points on several modular curves in \cite{krummthesis}.  One of these modular curves he studied was the hyperelliptic modular curve $Y_1(16)$. The hyperelliptic map  an infinite source of quadratic points on this modular curve. These quadraitc points are obtained as inverse images of rational points on $\P^1$. Krumm computed the class numbers of all quadratic fields occuring as the field of definition of a quadratic point correspong to a point in $\P^1$ of height at most 1000, with the restriction that the discriminant of this field did not exceed $10^{15}$. The data computed this way led him to ask the following question:
	
	\begin{question}[{\cite[2.6.13]{krummthesis}}] \label{q:krumm}
		Let $K \ncong \Q(\sqrt{-15})$ be an imaginary quadratic field over which $Y_1(16)$ has a
		point. Is it necessarily the case that the class number of $K$ is divisible by 10?
	\end{question}
	
	The main goal of this paper is to give an affirmative answer to the above question by proving the following:
	
	\begin{theorem}\label{thm:main}
		Let $K$ be an imaginary quadratic field such that $Y_1(16)(K) \neq \emptyset$  then the class number of $K$ is even. If furthermore $K \ncong \Q(\sqrt{-15})$ then the class number of $K$ is divisible by $10$.
	\end{theorem}
	
	The class number being divisible by $10$ is equivalent to it being divisible by $2$ and by $5$. The above theorem is proved by showing divisibility by $2$ and by $5$ separately.  
	
	Note that Krumm already gave the following equation
	\begin{align}
		\begin{split}
			X_1(16) : y^2 &= f_{16}(x) \\
			f_{16}(x) &:= x(x^2 + 1)(x^2 + 2 x - 1) \label{eq:x1_16}
		\end{split}
	\end{align}
	
	for $X_1(16)$. And $Y_1(16)$ is the open sub-variety  where $$x(x-1)(x+1)(x^2+1)(x^2-2 x -1)(x^2+2 x -1) \neq 0$$ and $x\neq \infty.$

	Krumm showed that all quadratic points on $Y_1(16)(\Q(\sqrt{d}))$ are of the form $(x,y) = (t,\pm\sqrt{f_{16}(t)})$ for some $t \in \Q$   \cite[Thm. 2.6.12]{krummthesis}. So more explicitly
	
	\begin{corollary}\label{cor:explicit}
		Let $t \in \Q$ such that $f_{16}(t) < 0$ and $t \neq -3, 1/3$ then $\Q(\sqrt{f_{16}(t)})$ is an imaginary quadratic field whose class number is divisible by $10$.
	\end{corollary}
	The values  $t=-3,1/3$ have to be excluded because these are precisely the values for which $\Q(\sqrt{f_{16}(t)})\cong\Q(\sqrt{-15})$, whose class number is $2$.
	
	The contrapositive of $\Cref{thm:main}$ shows that if $K \ncong \Q(\sqrt{-15})$ is an imaginary quadratic field whose class number is not divisible by 10 then $Y_1(16)(K) = \emptyset$. In \cite[\S 4]{banwait_derickx_quadratic_torsion} the squarefree integers $d$ with $|d| < 10,000$ such that $Y_1(16)(\Q(\sqrt{d})) \neq \emptyset$ were studied, and was successfully answered for all but 38 such integers $d$. This simple criterion on class numbers allows one to easily show that $Y_1(16)(\Q(\sqrt{d})) = \emptyset$ for 4955 out of the 6083 squarefree integers $d$ with $-10,000 < d <0$. In \cite[equation 4.1]{banwait_derickx_quadratic_torsion}  there is a list of the 38 values for which the methods in that paper fail.  The following corollary shows that the divisibility criterion is able to deal with 8 of those 38 unsolved cases.
	
	\begin{corollary}
		The class number of $\Q(\sqrt{d})$ is not divisible by $5$ and hence $Y_1(16)(\Q(\sqrt{d)}) = \emptyset$ for $$d \in \set{-7161, -6711, -6503, -6095,  -6005, -4847, -3503, -3199}.$$
	\end{corollary}
	
	\subsection{Divisibilty by 2.} The general approach for divisibility by $2$ is rather classical. Indeed genus field theory as developed by Gauss \cite{gauss66disquisitiones} implies that the narrow class number of a quadratic field is odd if and only if the absolute value of it's discriminant is a prime. So this case can be solved by analyzing the number of primes dividing the squarefree part of $f_{16}$ defined in \cref{eq:x1_16} below. This is done in \cref{sec:divisibility_by_2} by using two descent to relate the original question to the solution of certain Diophantine equations. As an auxiliary result the following theorem is proved: 

	\begin{theorem}\label{thm:diophantine_main}
	Let $f = a_4r^4+a_3r^3s+a_2r^2s^2+a_1rs^3+a_0s^4 \in \Z[r,s]$ be a squarefree homogenous polynomial with $a_0 \neq 0$,  then $y^2=f(x^2,z)$ has only finitely many solutions with $x,y,z\in \Z$ and $x,z$ coprime. 
	\end{theorem}
    
    In fact, \Cref{thm:diophantine_mainv2} is a generalization of this theorem to number fields. The proof  of the above Theorem is concrete enough that it allows for the explicit determination of the finitely many solutions in the example that is relevant for the divisibility by $2$ in case 2 in the proof of \Cref{prop:divisibilit_by_2}. The idea of the proof this theorem uses a from of descent. After reducing to the case $a_4=0$ one can show using descent that there is a finite set of integers $d_1,\ldots,d_n$ such that any solution to $y^2=f(x^2,z)$ is of the form $y^2=f(x^2,d_iw^2)$ for some $w$. So the solutions are covered by the rational points on finitely many genus $2$ curves.
	
	\subsection{Divisibility by 5.} The divisibility by $5$ is more difficult. The strategy can be outlined as follows.
	
	\begin{enumerate}
		\item Find an explicit  $[D] \in \Cl(X_1(16)_\Q)[5]  \cong \Pic(X_1(16)_\Q)[5]$ as a difference of cusps.
		\item Lift $[D]$ to  an element in $[\mathcal D] \in\Pic(X_1(16)_\Z)[5]$.
		\item Use the imaginary quadratic point $\mathcal P$ to construct the pull-back  $[\mathcal P^* (\mathcal D)]\in \Cl(\OO_{K})[5]$.
		\item Show that if $[\mathcal P^* (\mathcal D)]$ is trivial then $\mathcal P$ gives rise to a $\Q$-rational point on some auxiliary curve $Y_{16}$ over $\Q$ of genus $2$ and Mordell-Weil rank 3. This is done in \Cref{lem:auxiliary_curve}.
		\item Determine all rational poins in $Y_{16}(\Q)$ using a combination of two cover descent and elliptic curve Chabauty as described in \cite[\S 8]{bruin2009two}. This is done in \Cref{prop:auxiliary_points}.
	\end{enumerate}

The above strategy fits in a more general framework for transporting torsion elements in class groups of curves to class groups of number fields. This framework is worked out more generally in \Cref{sec:obstructions_to_divisibility}. In this more general setup there is a clear obstruction to carrying out step $(2)$. This obstuction can be described in terms of the minimal regular model of the curve, see \Cref{cor:lift_of_order_N} and \Cref{ex:non-principal}.  This obstruction happens to not be an issue for $X_1(16)$ because the chosen $[D]$ specializes to the connected component of the identity of the Neron model of $J_1(16)$.
The rational points in step $(4)$ can be interpreted as second obstruction coming from the fact that the order of a point might drop while transporting it form the class group of the curve to the class group of the number field. Note that steps $(1), (2)$ and $(3)$ can be carried out for more general number fields. However, step $(4)$ only works for imaginary quadratic fields. Indeed imaginary quadratic fields are the only fields for which one has that $\OO_K^*/\OO_K^{*5}={1}$. This fact is essential in interpreting obtruction in step $(4)$ in terms of $\Q$ rational points instead of $K$-rational points, see the usage of \Cref{lem:descent} and \Cref{prop:quadratic_to_rational} in \Cref{lem:auxiliary_curve}. Using the rational points $Y_{16}$ in step $4$ and $5$ once can show that the only cases in which $[\mathcal P^* (\mathcal D)]$ is trivial, are when the imaginary quadratic field over which $\mathcal P$ is defined is $\Q(\sqrt{-15})$ or $\Q(\sqrt{-2030})$. The reason that $\Q(\sqrt{-2030})$ doesn't show up as an exception in the main Theorem is because its class number is $40$, even though the explicit construction in step $(3)$ doesn't yield an element of order $5$.

 The idea for pulling back elements from the class group of the curve to the class group of a number field has also been used in \cite{GillibertLevenTorsiontoIdealClasses}.  However, in that article they do no pay any attention to the primes of bad reduction.  In order to get around the issues at the primes of bad reduction, they simply invert all those primes. Hence, applying their strategy to $X_1(16)$ one would end up with an element in $\Cl(\OO_{K}[\frac 1 2])[5]$ instead. The problem arises when one tries to show this element is nontrivial. Indeed, in general ${\OO_{K}[\frac 1 2]}^*/{\OO_{K}[\frac 1 2]}^{*5}$ could contain elements that are not coming from  ${\Z[\frac 1 2]}^*/{\Z[\frac 1 2]}^{*5}$ one runs into serious issues using an $\OO_{K}[\frac 1 2]$ version of \Cref{lem:descent} and \Cref{prop:quadratic_to_rational} to construct rational points on a curve defined over $\Q$ as in step $(4)$. And one will get points on some curve that depend on the quadratic field instead. And the fact that the curve from step $(4)$ doesn't depend on the field is essential to the argument.

\subsection{Structure of the paper}
In \Cref{sec:obstructions_to_divisibility} a general framework is outlined that allows one to study the obstructions of transporting elements in class groups of curves to class groups of number fields using rational points.

In \Cref{sec:divisibility_by_5} these techniques are applied to $X_1(16)$ to show the divisibility by $5$ part of the main Theorem.

In \Cref{sec:diophantine} the solution to a certain type of Diophantine problem is discussed that will pop up several times in the proof of the divisibilty by two part of the main theorem.

In \Cref{sec:divisibility_by_2} the divisibilty by $2$ part of the main Theorem is proved using a tedious 2 descent.

Finally, in \Cref{sec:real_quadratic_points} the parity of the narrow class numbers of real quadratic fields $K$ such that $Y_1(16)(K) \neq \emptyset $ is studied. Contrary to the imaginary quadratic case, there is a real quadratic field such that such that $Y_1(16)(K) \neq \emptyset$  with odd narrow class number whose discriminant has hundreds of digits. However, a heuristic argument shows that it is still expected that $Y_1(16)(K) \neq \emptyset$ for only finitely many real quadratic fields with odd class number.

\subsection{Acknowledgments}
I would like to thank Filip  for bringing this question to my attention and feedback on an earlier version of this paper. I would also like to thank David Krumm for the useful discussions and references to the literature.

\subsection{Data availability and reproducibility}

The computation for this article have been done using Magma (version V2.28-3)  \cite{magma}.
The code accompanying this paper has been released as part of the v0.2.1 release of mdmagma which is available at:
\begin{center}
	\url{https://github.com/koffie/mdmagma/tree/v0.2.1}
\end{center}
under the GNU General Public License v3.0 licence. 
All filenames given in the paper will refer to files in the directory \href{https://github.com/koffie/mdmagma/tree/v0.2.1/computations/derickx-X1_16_classgroups}{/computations/derickx-X1\_16\_classgroups} of this repository. The computations have been run on a server at the University of Zagreb with an Intel Xeon W-2133 CPU @ 3.60GHz with 12 cores and 64GB of RAM running Ubuntu 18.04.6 LTS. 

The computations for the divisibility by 2 of the class number took around 5 seconds and used less then  250 MB of RAM.
The computations for the divisibility by 5 of the class number took less than 5 minutes and used less then 300 MB of RAM.

Logs containing the output of the computations are avaiable in the \githubbare{logs} subfolder. More information on how to reproduce these results can be found in the \githubbare{README.md} file.

Claims in this article that are based on computations have been labeled in the source code in order to make it easy to navigate between this text and the code.
 For example the proof of \cref{lem:auxiliary_curve} starts with the sentence: "Let $P_1 := (1,2)$, then  $P_1 - \omega(P_1)$ is of order $5$ \github{order_5.m}{Claim 1}."
 
 In order to verify this claim one can open the file\githubbare{order_5.m} and search for "Claim 1" to find the following code that supports this claim.
 \begin{center}
 \begin{verbatim}
// Claim 1: p1 - w(p1) is or order 5, note p2 = w(p1);
is_of_order_5, g := IsPrincipal(5*(p1-p2));
assert is_of_order_5;
print "Claim 1 successfully verified";
\end{verbatim}
 \end{center}

The results of the above claim can also easily be seen in the corresponding log file \githubbare{logs/order_5.txt}.

	\section{Obstructions to Divisibility of class groups}\label{sec:obstructions_to_divisibility}
	Throughout this section $X$ is a nice\footnote{nice meaning smooth projective and geometrically integral} curve of genus $\geq 2$ over $\Q$ and $\mathcal X$ is it's minimal regular model over $\Z$.  Additionally $K$ will be a number field whose maximal order we denote by $\OO_K$.
	
	According to \cite[Prop. 7.2.16]{liu_geom_curves} the notions of Cartier and Weil divisors on $\mathcal X$ agree, the same also holds for $X$.  Their class groups will be denoted by $\Cl(\mathcal X)$ and $\Cl(X)$ and their divisor groups by $\Div(\mathcal X)$ and $\Div(X)$. Note that $\Q(X) \cong \Q(\mathcal X)$. If $g \in \Q(\mathcal X)$ then $\ddiv(g)$ will denote the associated principal divisor on $\mathcal X$, and $\ddiv(g)_\Q$ the one on $X$. For every divisor $\mathcal D$ respectively $D$ on $\mathcal X$ resp. $X$ there  are the associated line bundles $\OO_{\mathcal X} (\mathcal D)$ and $\OO_X(D)$. Sending a divisor to its associated line bundle induces isomorphisms $\Cl(\mathcal X) \to \Pic \mathcal X$ and $\Cl(X) \to \Pic X$ \cite[Corollary 7.1.19]{liu_geom_curves}.
	
	Let $P \in X(K)$ be a point, then this point uniquely extends to a point $\mathcal P \in \mathcal X(\OO_K)$.
	
	The group structures of $\Cl(X)$ and $\Cl(\mathcal X)$ are written additively, however the group structure of $\Cl(\OO_K)$ is written multiplicatively. The reason for this is inconsistency is that $\Cl(X)$ and $\Cl(\mathcal X)$ are usually  thought of as divisors modulo principal divisors,  which are usually denoted additively. While $\Cl(\OO_K)$ are ideals modulo principal ideals, which are usually denoted multiplicatively. This leads to equations like $\mathcal P^*(\mathcal D)^N=\mathcal P^*(N\mathcal D)$ and $\mathcal P^*(0)=1$ rather than $N\mathcal P^*(\mathcal D)=\mathcal P^*(N\mathcal D)$ and $\mathcal P^*(0)=0$.

	Let $\mu: X \to \mathcal X$ be the inclusion then pullback defines group homomorphisms $\mu^* : \Div(\mathcal X) \to \Div(X)$ and $\mu^* : \Cl(\mathcal X) \to \Cl(X)$. Sometimes $\mathcal D_\Q$ will also be used as shorthand notation for $\mu^*(\mathcal D)$. The kernel of $\mu^* : \Div(\mathcal X) \to \Div(X)$ exactly consists of the vertical divisors, i.e. the divisors supported on the closed fibres of $\mathcal X \to \Spec \Z$ . If $D \in \Div X$, then $D$ can be written as the difference of two effective divisors $D_1 - D_2$, we will use $\overline D$ as a shorthand notation for $\overline{D_1}- \overline{D_2}$ where $\overline{D_i}$ denotes the schematic closure of $D_i$ in $\mathcal X$. The horizontal divisors on $\mathcal X$ are exactly the divisors of the form $\overline D$. Note that $\mu^*(\overline D) = D$ so that $\mu^*$ is surjective. While $D \to \overline D$ defines a homomorphism $\Div(X) \to \Div(\mathcal X)$, it doesn't induce a homomorphism $\Cl(X) \to \Cl(\mathcal X)$ since $D \to \overline D$ doesn't necessarily preserve principal divisors, see \Cref{ex:non-principal} for an example.
	
	Consider the diagram: $$\xymatrix{
		&\Cl(\mathcal X) \ar[dl]_{\mathcal P^*}\ar[dr]^{\mu^*}& \\
		\Cl(\OO_K) &  & \Cl(X)
	}$$
	Let $N>0$ be an integer and $[D] \in \Cl(X) \cong \Pic X$ be the linear equivalence class associate to a divisor $D$ and assume that the order of $[D]$ is exactly $N$. The main question we address in this section is when can we use $\mu^*$ and $\mathcal P^*$ to construct an element of $\Cl(\OO_K)$ of order $N$.

	There are two obstructions to constructing this element of order $N$ in $\Cl(\OO_K)$.  The first one is that there might not be a $[\mathcal D] \in \Cl(\mathcal X)$ of order $N$ such that $\mu^*([\mathcal D]) = [D]$. Then second problem is that $\mathcal P*$ might not be injective.
	
	\Cref{lem:extending_divisor_classes} and it's corollary give sufficient conditions for conditions for $[D]$ to have a lift of order $N$. However, before the statement of this lemma, we first need to formulate a condition that allows us to deal with the bad fibres.
	
	\begin{definition}
	Let $p$ be a prime and $\Gamma_1, \ldots, \Gamma_r$ be the irreducible components of 	$\mathcal X_{\F_p}$. Let $d_1, \ldots, d_r$ denote the multiplicities of   $\Gamma_1, \ldots, \Gamma_r$ in $\mathcal X_{\F_p}$. If $gcd(d_1, \ldots, d_r)=1$ then $\mathcal X$ is said to have {\it coprime multiplicities at the fibre above $p$}. If $\mathcal X$ has  { coprime multiplicities at fibre above $p$} for all primes $p$, then $\mathcal X$ it is said that {\it the fibres of $\mathcal X$ have  coprime multiplicities}
	\end{definition}

\begin{lemma}\label{lem:coprime_multiplcities}
	Let $p$ be a prime such that $\mathcal X$ has coprime multiplicities at the fibre above $p$, and let $D$ be a divisor on $\mathcal X$ supported at the fibre above $p$. If $D \cdot C = 0$ for all divisors $C$ in the fibre above $p$, then there is an integer $m$ such that $D=m \mathcal X_{\F_p}$. 
\end{lemma}
\begin{proof}
Let $\Gamma_1, \ldots, \Gamma_r$ be the irreducible components of $\mathcal X_{\F_p}$. Let $d_1,\ldots,d_r,\allowbreak a_1, \ldots, a_r \in \Z$ such that $\mathcal X_{\F_p} = \sum_{i=1}^r d_i \Gamma_i$ and $D = \sum_{i=1}^r a_i \Gamma_i$ . By \cite{liu_geom_curves}[Thm 9.1.23] there exists an $m \in \R$ such that $D =m \mathcal X_{\F_p}$, implying that for $1 \leq i \leq r$ one has $a_i = \gamma d_i$. Since the $d_i$ and $a_i$ are integers, in fact $m\in \Q$. Since $\gcd(d_1,\ldots, d_r) =1$ it follows that $m \in \Z$. 
\end{proof}

When $p$ is a prime of good reduction for $\mathcal X$, then $\mathcal X_{\F_p}$ consists of a single component of multiplicity 1. In particular, $\mathcal X$ has coprime multiplicities at the fibre above  primes of good reduction. So the primes of bad reduction are the only primes at which the fibres of $\mathcal X$  can fail to have coprime multiplicities.

	\begin{lemma}\label{lem:extending_divisor_classes}
		Let $N$ be an integer and $D$ be divisor on $X$ such that $[D] \in Cl(X)$ is of order $N$. Let $g \in \Q(\mathcal X) =\Q(X)$ such that $ND = \ddiv(g)_\Q$. Then there exist a $\gamma \in \Q^*$ and divisor $\mathcal D' \in \Div(\mathcal X)$ supported at the fibres of bad reduction  $N\overline D + \mathcal D' = \ddiv(\gamma g)$. If furthermore the fibres of $\mathcal X$ have  coprime multiplicities and $\overline D \cdot C = 0$ for every vertical divisor $C$ on $\mathcal X$ then $\gamma$ can be chosen such that $\mathcal D'  = 0$.
	\end{lemma}
	\begin{proof}
	Since $[D]$ is of order $N$ one take $g\in \Q(\mathcal X)$  to such that $ND = \ddiv(g)_\Q$. Let $ \mathcal D_0 := \ddiv(g)-N\overline D $. Since $\mathcal D_{0,\Q} =  \ddiv(g)_\Q-ND=0$ one has that $\mathcal D_0$ is vertical, although it might still contain fibres of good reduction in it's support. For every prime $p$ of good, let $m_p$ denote the mulitplicity of  $\mathcal X_{\F_p}$ in  $\mathcal D_0$. Set $m_p=0$ for the primes of bad reduction.  Then by defining $\gamma := \prod_p p^{-m_p}  $  and taking $\mathcal D' := \ddiv( \gamma g)-N\overline D = \mathcal D_0 + \ddiv(\gamma)$ one can ensure that $\mathcal D'$ has no fibres of good reduction in it's support.
	
	For the second part, assume $\overline D \cdot C = 0$ for every vertical divisor $C$ on $\mathcal X$. Then $\mathcal D_0 \cdot C = \ddiv(g)\cdot C-N\overline D\cdot C = 0$, because $ \ddiv(g)\cdot C = 0$ by \cite{liu_geom_curves}[Thm. 9.1.12 (c)]. In particular if one defines $\mathcal D_{0,p}$ to be the part of $\mathcal D_0$ supported at the fibre above  a prime $p$, then \Cref{lem:coprime_multiplcities} ensures the existence of an integer $m$ such that $\mathcal D_{0,p} = m \mathcal X_{\F_p}$. For every prime $p$, let $m_p$ be the integer such that $\mathcal D_{0,p} = m_p \mathcal X_{\F_p}$. Then by setting $\gamma := \prod_p p^{-m_p}$ as before one can assure that $\mathcal D':= \ddiv( \gamma g)-N\overline D = \mathcal D_0 + \ddiv(\gamma)$ has no vertical fibres in it's support. However, since $\mathcal D'$ is vertical this  implies $\mathcal D'=0$.

	\end{proof}
	
	\begin{corollary}\label{cor:lift_of_order_N} If $[D]$ is of order $N$ and $\overline D \cdot C = 0$ for every fibral divisor $C$ then $\overline D$ is a divisor on $\mathcal X$ that maps to $D$ under $\mu^*$ and whose class $[\overline D]$ is of order $N$.
	\end{corollary}
	The condition that $\overline D \cdot C = 0$ for every fibral divisor $C$ cannot be removed from the corollary as is shown by the example below.
	\begin{example}\label{ex:non-principal}
		If $p \equiv 1 \mod 12$ is a prime and $p \neq 13$ then the special fibre of $X_0(p)$ over $\F_p$ consists of two components $C_0$ and $C_1$ intersecting each other in $\frac {p-1} {12}$ points \cite[Diagram 1 p. 63]{mazur_mod_curves}. $X_0(p)(\Q)$ contains two cusps $0_\Q$ and $\infty_\Q$, and their difference $D := 0_\Q - \infty_\Q$ is of order $N := \frac {p-1} {12}$ in $\Cl(X_0(p)_\Q)$. Their specializations $0_{\F_p}$ and $\infty_{\F_p}$ lie in different components. Let's assume we numbered them such that $0_{\F_p}$ lies on $C_0$ and $\infty$ lies on $C_1$.  Let $g \in \Q(X_0(p))$ such that $D= \ddiv(g)_\Q$. Then we have the following equalities of intersection multiplicities $\ddiv(\gamma g) \cdot C_i =0$, $C_i \cdot C_i = -N$, $C_0 \cdot C_1 = N$, $C_0 \cdot \overline D = 1$  and $C_1\cdot  \overline D = -1$. It follows that any divisor $\mathcal D'$ of the form $\mathcal D' = \ddiv(\gamma g)-N\overline D$ has to satisfy $\mathcal D' \cdot C_0 = -N$ and $\mathcal D' \cdot C_1 = N$. In particular there is an integer $k$ such that $\mathcal D'$ is of the form $C_0 + k(C_0 + C_1) = C_0 + \ddiv(p^k)$. Since the order of $[C_0] = [C_0 + k(C_0 + C_1) ]$ in $\Cl(X_0(p))$ is infinite, it is impossible to add a vertical divisor to $\overline D$ in order to make it of order $N$. In particular, there is no divisor $E$ on $\mathcal X$ such that $\mu^*(E) = D$ and $[E]$ is of order $N$. These difficulties occur because file $ND$ is a principal divisor, but $N\overline D$ is not.
	\end{example}
	
	When $\mathcal D'$ is nonzero one can still use \Cref{lem:extending_divisor_classes} to show that for certain points $P$ that $\mathcal P^*(\overline{D})$ is of order dividing $N$.  This can be done when $\mathcal P^*(\mathcal D')$ is principal. One way in which this happens is if the image of $\mathcal P$ doesn't meet the support of $\mathcal D'$ in which case $\mathcal P^*(\mathcal D')=1$.
	
	\begin{proposition}\label{prop:descent}
		Let $N>0$ be an integer and $\mathcal D, \mathcal D'$ be divisors on $\mathcal X$ such that $N\mathcal D +\mathcal D' = \ddiv(g)$ is a principal divisor.  Let  $\mathcal P \in \mathcal X(\mathcal O_K)$ be a point not contained in the support of $\mathcal D$ and $\mathcal D'$ such that $\mathcal P^*(\mathcal D') = 1$. Then the order of $[\mathcal P^*(\mathcal D)]$ divides $N$. Denote the order of $[\mathcal P^*(\mathcal D)]$ by $N'$, then there are $\alpha \in K^*$ and $\beta \in \OO_K^*$ such that $\beta g(\mathcal P)=\alpha^{N/N'}$.
	\end{proposition}
	
	\begin{proof}
		The order of $[\mathcal P^*(\mathcal D)]$ dividing $N$ follows from the computation below
		\begin{align*}
		\mathcal P^*(\mathcal D)^N= \mathcal P^*(\mathcal D)^N+ \mathcal P^*(\mathcal D')=\mathcal P^*(N\mathcal D + \mathcal D') = \mathcal P^*(\ddiv(g)) = \ddiv(g(\mathcal P))
		\end{align*}
		showing that $\mathcal P^*(\mathcal D)^N$ is a principal fractional ideal.
		
		For the second part of the statement, let $\alpha$ be a generator of the principal fractional ideal $ \mathcal P^*(\mathcal D)^{N'}$. Then $\alpha^{N/N'}$ and $g(\mathcal P)$ are both generators of the ideal $\mathcal P^*(\mathcal D)^N$, implying that $\alpha^{N/N'}/g(\mathcal P)$ is an element of $\OO_K^*$. The second part follows by setting $\beta := \alpha^{N/N'}/g(\mathcal P)$.
	\end{proof}
	
	\subsection{Relation with descent}\label{sec:descent}
	
	Let $g \in \Q(X)^*$, $\beta \in K^*$ and $d$ an integer. Then $\Q(X)[z]/(z^d-\beta g) = \prod_{i}^n L_i$ is a product of fields of transcendence degree $1$ over $\Q$. So each $L_i$ is the function field of a smooth and projective geometrically reduced curve $Y_i$. These curves $Y_i$ are not necessarily geometrically irreducible. Since each $L_i$ contains $\Q(X)$, the curves $Y_i$ come equipped with a map $Y_i \to X$. We define the curve $Y_{d,g,\beta}$ as follows:
	
	\begin{align}Y_{d,g,\beta} := \coprod_{i=1}^n Y_i
	\end{align}
	
	If $\beta,\beta' \in K^*$ such that $\beta/\beta' \in K^{*d}$ then $Y_{d,g,\beta}$ and $Y_{d,g,\beta'}$ are isomorphic via an isomorphism that respects the map to $X$.  So that the isomorphism class of $Y_{d,g,\beta}$ only depends on equivalence class of $\beta$ in  $K^*/K^{*d}$.
	
	\begin{lemma}\label{lem:descent}
		Let $N>0$ be an integer and $\mathcal D, \mathcal D'$ be divisors on $\mathcal X$ such that $N\mathcal D +\mathcal D' = \ddiv(g)$ is a principal divisor.  Let  $\mathcal U := \mathcal X \setminus \Support(\mathcal D')$ and $\mathcal P \in \mathcal U(\OO_K)$ be a point not contained in the support of $\mathcal D$. Let $M = gcd(N, \#\Cl(\OO_K))$ and $S \subset \OO_K^*$ be a set of coset representatives of $\OO_K^*/\OO_K^{*N/M}$. Then $\mathcal P_K$ is in the image of $$\coprod_{\beta \in S} Y_{N/M,g,\beta}(K) \to X(K).$$
	\end{lemma}
	
	The power of the above lemma really depends on $\Support(\mathcal D')$. Ideally one would try to arrange for $\mathcal D' =0$ so that $\Support(\mathcal D') =0$ and hence $\mathcal U(\OO_K) = \mathcal X(\OO_K) $ so  the lemma is applicable to almost all points of $X(K)$. On the other extreme, if there is a prime $p$ such that $\mathcal X_{\F_p} \subset \Support(\mathcal D')$, then $\mathcal U(\OO_K) = \emptyset$ and the lemma has no content.
	\begin{proof}
		By definition of $\mathcal U$ the image of $\mathcal P$ doesn't meet the the support of $\mathcal D'$, so that $\mathcal P^*(\mathcal D')=1$. In particular the conditions of \Cref{prop:descent} are satisfied. Let $N'$ be the order of $\mathcal P^*(\mathcal D)$, then $N' \mid M$. The proposition gives us $\beta' \in \OO_K^*$ and $\alpha \in K^*$ such that $\beta' g(\mathcal P)=\alpha^{N/N'}$. Let $\beta \in S$ such that $\beta=\beta' \gamma^{N/M}$ for some $\gamma \in \OO_K^*$. Then setting $z=\gamma\alpha^{M/N'}$ one has $y^{M/N} =\gamma^{N/M}\alpha^{N/N'}=\gamma^{N/M}\beta' g(\mathcal P) =\beta g(\mathcal P) .$ So taking $z=\gamma\alpha^{M/N'}$ gives us a point in $Y_{N/M,g,\beta}(K)$ mapping to $\mathcal P_K$.
	\end{proof}

	The above lemma can be interpreted as a form of descent. Indeed there is an action of $\mu_{N/M}$ acts on $Y_{N/M,g,1}$ by sending $z$ to $\zeta_{N/M}z$. This action is such that $X =  Y_{N/M,g,1}/\mu_{N}$. From the covering $Y_{N/M,g,1} \to X$ one can construct a $\mu_{N/M}$-torsor $\mathcal Y \to \mathcal U$. Since $S \subset  K^*/K^{*N/M} =H^1(\Gal(K),\mu_{N/M}(\overline K))$ one can view $\coprod_{\beta \in S} Y_{N/M,g,\beta}$ as a union of twists of $Y_{N/M,g,1}$ over $X$. And the above lemma says that $\mathcal P_K$ has to come from one of the twists in this union.
	
	\begin{remark}\label{rem:imaginary_quadratic}
		The above lemma is more general than what we need. Indeed for $X_1(16)$ we will be able to arrange it in such a way that $\mathcal D' = 0$ so that $\mathcal U = \mathcal X$, and hence that the only restriction on $\mathcal P \in \mathcal X(\OO_X) = X(K)$ is that it shouldn't be one of the finitely many points contained in the support of $\mathcal D$. Additionally, when restricted to $N$ odd and $K \ncong \Q(\zeta_{3})$ an imaginary quadratic field then $\OO_K/\OO_K^{*N/M} = \set{1}$ so that one can take $S = \set{1}$ and hence $\mathcal P_K$ lifts to a point on $Y_{N/M,g,1}$.
	\end{remark}
	
	\subsection{Imaginary quadratic points on hyperelliptic curves}
	
	This section contains a construction that specializes the results obtained earlier to a rather special situation. However this situation is exactly what happens for $X_1(16)$, and helps explaining the divisibility by $5$ part of the phenomenon that motivated Krumm to ask \Cref{q:krumm}.
	
	Let $f \in \Q[x]$ of degree $>4$, then we use $X_f$ to denote the hyperelliptic curve given by
	\begin{align}
	y^2 = f(x)
	\end{align}
	and $\mathcal X_f$ it's minimal regular model. The hyperelliptic involution of $X_f$ will be denoted by $\omega$.
	
	\begin{proposition}\label{prop:extending_divisor_classes_hyperelliptic}
		Let $N>0$ be an integer and $D$ be a divisor on $X_f$ such that $[D -\omega(D)] \in \Cl(X_f)$ is of order $N$. Assume that $\overline D \cdot C = \overline{\omega(D)} \cdot C$ for every vertical divisor $C$. Then there exists an $g(x,y) \in \Q(\mathcal X)$ such that $N(\overline{D-\omega(D)}) = \ddiv g$. Additionally there is an $s \in \set{-1,1}$ such that $g(x,-y) = g \circ \omega = sg(x,-y)^{-1}$.
	\end{proposition}
	
	\begin{proof}
		The existence of $g$ is \Cref{lem:extending_divisor_classes}. Now let  $s := g(x,y)/g(x,-y)$, it suffices to show $s = \pm 1$. Since $\ddiv s =\ddiv g(x,y)/g(x,-y)  = N(\overline{D-\omega(D)}) - \omega(N(\overline{D-\omega(D)}) )=0$. This means that $s\in \Q$ is a constant function.  However the support of $\ddiv s$ also doesn't contain any vertical divisors, so $s$ has to be a unit  in $\Z$ so that $s = \pm 1$.
		
	\end{proof}
	
	Based on this, let $s=-1$ or $1$ and suppose that $d$ is odd when $s=-1$. Let $g(x,y) \in \Q(X_f)$ such that $g(x,-y)=sg(x,y)^{-1}$. Then the hyperelliptic involution of $X_f$ lifts to an involution of $Y_{d,g,1}$ in the following way:
	\begin{align}
	\begin{split}
	\tau : \quad\, Y_{d,g,1} &\to Y_{d,g,1} \\
	(x,y,z) &\mapsto (x,-y,sz^{-1})
	\end{split}
	\end{align}
	
	In the light of \Cref{rem:imaginary_quadratic}, in order to study quadratic imaginary quadratic points on $X_f$ where the class group is not divisible by $N$, it suffices to know the quadratic points on $Y_{d,g,1}$ for the divisors $d|N$ that are strictly larger then  $1$. Now often, most quadratic points come from the hyperelliptic involution and hence have an $x$-coordinate that is in $\Q$. The following proposition helps simplify the above problem by turning it into a question about $\Q$-rational points on a curve instead of $K$-rational points.
	
	\begin{proposition}\label{prop:quadratic_to_rational}
		Let $s=\pm -1$ and $d \neq 1$ be an odd integer. Let $K$ be a quadratic field such that $\zeta_d^i \notin K$ for $1 \leq i <d$. Let $g(x,y) \in \Q(X_f)$ such that $g(x,-y)=s g(x,y)^{-1}$  and $P \in Y_{d,g,1}(K)$ be a quadratic point such that $x(P)\in \Q$ and $g(P) \neq 0,\infty$. Then $P$ maps to a $\Q$-rational point on either $X_f$ or $Y_{d,g,1}/\tau$.
	\end{proposition}
	The condition $g(P) \neq 0,\infty$ ensures that there is no need to worry about the desingularisation that is happening in the definition of $Y_{d,g,1}$. Indeed, the singularities of the covering of $X_f$ given by $z^d=g(x,y)$ only happen at the poles and zero's of $g$.
	\begin{proof}
		Let $\iota : K \to K$ denote the involution of $K$ and denote  the image of $P$ under $Y_{d,g,1} \to X_f$ by $Q$.  Let $(\alpha,\beta,\gamma)$ be the coordinates of $P$ so that $Q=(\alpha,\beta)$ and $\alpha \in \Q$. If $Q$ is a $\Q$-rational point there is nothing to proof. So we assume $Q$ is not $\Q$-rational. This means that $K=\Q(\sqrt{f(\alpha)})$,  $\beta = \pm \sqrt{f(\alpha)}$ and $\iota(Q) = (\alpha, -\beta)$. Now let $\gamma'$ be such that $\iota(P)=(\alpha, -\beta,\gamma')$ and $\gamma'' := s\gamma^{-1}$ be such that $\tau(P)=(\alpha, -\beta,\gamma'')$. This means that $\iota(P)$ and $ \tau(P)$ both map to $\iota(Q)=(\alpha,-\beta)$ on $X_f$. In particular $(\gamma''/\gamma')^d=g(\alpha,-\beta)/g(\alpha,-\beta)=1$ so that $\gamma''/\gamma'=\zeta_d^i$ for some $0 \leq i < d$. Since $\zeta_d^i \not \in K$ when $i >0$ we have $\zeta_d^i=1$ so that $\gamma''=\gamma'$ and actually $\tau(P) = \iota(P)$. In particular the image of $P$ in $Y_{d,g,1}/\tau$ is invariant under the action of $\Gal(K/\Q)$ and hence $\Q$-rational.
	\end{proof}

	\section{Divisibility by 5}
	\label{sec:divisibility_by_5}
	In this section $\mathcal X_1(16)$ will denote the minimal regular model of $X_1(16)$ over $\Z$.
	
	\begin{theorem}\label{thm:divisibility_by_5}
		Let $K\ncong \Q(\sqrt{-15})$ be an imaginary quadratic number field. If $Y_1(16)(K) \neq \emptyset$, then the class number of $K$ is divisible by $5$.
	\end{theorem}
	The first step in the proof of this theorem is to show that the exceptions give rational points on an auxiliary hyperelliptic curve.
	
Define the curve $Y_{16}$ and the function $h_{16} \in \Q(Y_{16})$ by 
\begin{align}
Y_{16} : y^2=& z^6 + z^5 - 5 z^3 + z + 1 \label{eq:y16} \\
h_{16}(z,y) :=& \frac {(2z^2 - 6 z + 2)y + (10 z^2 - 10 z - 
	2)}{(z^5 - 5 z^4 + 5 z^3 + 5 z^2 - 5 z - 3)} -1 \label{eq:h16}
\end{align}
	\begin{lemma}\label{lem:auxiliary_curve}

Let $K$ be an imaginary quadratic number field whose class number is not divisible by 5 and $Q \in Y_1(16)(K)$ then there is a  point $R \in Y_{16}(\Q)$ with $h_{16}(R) = x(Q)$.
	\end{lemma}

	\begin{proof}
		Let $P_1 := (1,2)$, then  $P_1 - \omega(P_1)$ is of order $5$ \github{order_5.m}{Claim 1}. The prime 2 is the only prime of bad reduction. So $\overline {P_1}  \cdot C = 1 = \overline{\omega(P_1)} \cdot C$ for every irreducible vertical divisor lying above a prime $>2$. The fibre $\mathcal X_1(16)_{\F_2}$ consists of two irreducible components \github{order_5.m}{Claim 3} which will be denoted by $C_1$ and $C_2$. The points $\overline {P_1} $ and $\overline{\omega(P_1)}$ reduce to the same component modulo $2$ \github{order_5.m}{Claim 4}, by relabeling we can assume this component is $C_1$  so that  $\overline {P_1}  \cdot C_1 = 1 = \overline{\omega(P_1)} \cdot C_1$ and  $\overline {P_1}  \cdot C_2 = 0 = \overline{\omega(P_1)} \cdot C_2$. In particular the conditions of \Cref{prop:extending_divisor_classes_hyperelliptic} are satisfied for $D=P_1$ and $N=5$. The function
		\begin{align*}
		g(x,y)  := \frac {(-6x^2 - 4 x + 2)y + (x^5 + 13 x^4 - 2 x^3 + 10 x^2 - 7 x + 1)}{x^5 - 5 x^4 + 10 x^3 - 10 x^2 + 5 x - 1}
		\end{align*}
		is such that $\overline {P_1}  - \overline{\omega(P_1)} = \ddiv g(x,y)$ \github{order_5.m}{Claim 2}.
		For an imaginary quadratic field $\OO_K^*$ is isomorphic to $\Z/n\Z$ with $n \in {2,4,6}$. In particular, $\OO_K^*/\OO_K^{*5}=1$. Applying \Cref{lem:descent} with $N=5$, $M=1$ and $S=\set{1}$ yields that $Q$ lifts to a quadratic point $P \in Y_{5,g,1}(K)$.
		
		By \cite[Thm. 2.6.12]{krummthesis} the point $Q \in Y_{16}(K)$ has to have $x(Q) \in \Q$, so that $x(P)=x(Q) \in \Q$. In particular, the conditions of \cref{prop:quadratic_to_rational} are satisfied. So we get that $P$ maps to a rational point on $Y_{5,g,1}/\tau$. One can compute that $Y_{5,g,1}/\tau \cong Y_{16}$ \github{order_5.m}{Claim 5}. Since $x$ is a function on $Y_{5,g,1}$ that is invariant under $\tau$ it factors via a unique function $h_{16}$ on $Y_{16}$. This function $h_{16}$ is given by \cref{eq:h16}  \github{order_5.m}{Claim 6}.
	\end{proof}
	

	\begin{proposition}\label{prop:auxiliary_points}
		The rational points on $Y_{16}$ (\cref{eq:y16}) are the $6$ points
		$$(0, -1), (0, 1), (1/3, -29/27),
		(1/3, 29/27), (3, -29), (3, 29) $$
		together with the two points at $\infty$.
	\end{proposition}
   The Mordell-Weil rank of $J(Y_{16})(\Q)$ has rank 3 \github{order_5.m}{Claim 7}. So that just using normal Chabauty or even quadratic Chabauty is not possible to for determining all rational points on this curve. The automorphism group of of $Y_{16}$ is $(\Z/2\Z)^2$ \github{order_5.m}{Claim 8}, making it bi-elliptic. It has a map of degree 2 to the elliptic curve 
   \begin{align*}
   E: y^2 = x^3 - x^2 + 17 x - 13
   \end{align*}
   which has rank 2 \github{order_5.m}{Claim 9 and 10}. It follows that $J(Y_{16})$ decomposes up to isogeny into two elliptic curves. One of rank 1, and another of rank 2. So these individual factors can also not be used to determine $Y_{16}(\Q)$. 
	
	\begin{proof} This was done by combining two cover descent \cite{bruin2009two} with elliptic curve Chabauty \cite{bruin2003chabauty}. For more details on how to combine these two methods see \cite[\S 8]{bruin2009two} or \cite[\S 4.2]{banwait_derickx_quadratic_torsion}. 
		What follows here is a short summary on how the computations were carried out. Let $\alpha$ be a root of $x^3 - x^2 + 2 x + 2$ and define $K$ to be the cubic number field $\Q(\alpha)$. By \github{order_5.m}{Claim 11} the defining hyperelliptic polynomial for $Y_{16}$ in \cref{eq:y16} factors as 
	\begin{align*}
	 z^6 + z^5 - 5 z^3 + z + 1 =& g_1 g_2 \\
	 g_1 :=& z^2 + (-\alpha^2 + \alpha - 1)z + 1\\
	 g_2 := & z^4 + (\alpha^2 - \alpha + 2)(z^3+z) + (\alpha^2 - 3\alpha + 3)z^2  + 1
	\end{align*}
	
	For $\gamma \in K$ let $E_{\gamma}$ be the genus 1 curve given by $y^2=\gamma g_2$. The outcome of the two cover descent is that every rational point on $Y_{16}$ either shares the same $z$ coordinate with a point on $E_{1}$ or with a point on $E_{\alpha^2+2}$. Applying elliptic curve Chabauty to the curves $E_{1}$ and $E_{\alpha^2+2}$ shows that the rational points on $Y_{16}$ are indeed as claimed in the proposition \github{order_5.m}{Claim 12}. The two points $(0,-1)$ and $(0,1)$ together with the two points at infinity come from $E_{1}$ the other four come from $E_{\alpha^2+2}$.
		
	\end{proof}

\begin{proof}[Proof of \Cref{thm:divisibility_by_5}]
   \Cref{lem:auxiliary_curve} reduces the proof of \Cref{thm:divisibility_by_5} to the computation of all points on $Y_{16}$. These rational points were determined in \Cref{prop:auxiliary_points}. These $8$ points give 8 different $x$ coordinates on $Y_1(16)$ that we need to check for potential counter examples to \Cref{thm:divisibility_by_5}. The corresponding points on $Y_1(16)$ together with the fields and class numbers are computed in \github{order_5.m}{Claim 13} and listed in \cref{tab:point_info}. From that table it follows that $\Q(\sqrt{-15})$ is indeed the only imaginary quadratic field where the class number is not divisible by $5$.
   \end{proof}

 The other imaginary quadratic field $\Q(\sqrt{-2030})$ listed in \cref{tab:point_info}  does have a class number 5 that is divisible by 5. However, at the moment I have no explanation for this divisibility by $5$. Indeed, the constructed element $\mathcal P^*(\mathcal D) \in \Cl(\OO_K)$ in the proof of \Cref{lem:descent} is trivial for all $\mathcal P \in \mathcal X_1(16)(\Q(\sqrt{-2030})))$ \github{order_5.m}{Claim 14}. So the divisibility by $5$ of $\Cl(\Q(\sqrt{-2030}))$ is not explained by the construction in $\cref{sec:obstructions_to_divisibility}$.
\begin{table}
	\centering

   \begin{tabular}{|r|r|r|r||}
    \hline 
   	Point on $Y_{16}$ & Points on $Y_1(16)$ & K & $\#\Cl(K)$ \\ 
   	\hline
   	    (1 : 1 : 0) & $(1,\pm 2) $ & $\Q$ & 1\\ 
        (0 : 1 : 1) & $(-1,\pm 2)$ & $\Q$ & 1\\ 
   	\hline 
   	   	(3 : 29 : 1) &$ \infty$ & $\Q$ & 1\\ 
   	   	(1 : 29 : 3) &  (0, 0) & $\Q$ & 1 \\
   	\hline 
   	    (1 : -1 : 0) & $(-3, \pm 2\sqrt{-15})$ & $\Q(\sqrt{-15})$ & 2\\ 
   	    (0 : -1 : 1) & $(1/3, \pm {2}/{27}\sqrt{-15 })$ & $\Q(\sqrt{ -15})$ & 2\\ 
   	\hline 
        (3 : -29 : 1) & $(-242/29, \pm 94721/24389\sqrt{-2030 })$ & $\Q(\sqrt{ -2030})$ & 40\\ 
   	    (1 : -29 : 3) & $(29/242, \pm8611/ 1288408\sqrt{-2030 })$ & $\Q(\sqrt{-2030}) $ & 40\\ 
   	\hline 
   \end{tabular} 

   \caption{Points where divisibility by $5$ could potentially fail. \label{tab:point_info}}
\end{table}
   
\section{A Diophantine problem}\label{sec:diophantine}

The following Theorem implicitly shows up in the divisibility by two part. That is, in case 2 in the proof of \Cref{prop:divisibilit_by_2} we need to solve a Diophantine problem of this shape. The goal of this section is to show more generally that such Diophantine problems only have finitely many solutions.

\begin{theorem}\label{thm:diophantine_mainv2}
	Let $K$ be a number field and let $f = a_4r^4+a_3r^3s+a_2r^2s^2+a_1rs^3+a_0s^4 \in \OO_K[r,s]$ be a squarefree homogeneous polynomial with $a_0 \neq 0$,  then $y^2=f(x^2,z)$ has only finitely many solutions $(x:y:z)\in \PP^2(1,4,2)(\OO_K)$ in weighted projective space with $x,z$ coprime. 
\end{theorem}
A more elementary way of phrasing this is that up to scaling transformations of the form $(x,y,z)\mapsto (ux:uy^4,u^2z)$ there are only finitely many solutions with $x,y,z \in \OO_K$ and $x,z$ coprime.

\begin{proof}
	After replacing $K$ an extension of degree at most $4$ we may and do assume that $f$ has a linear factor $b_1r+b_0s$. Since $a_0 \neq 0$ also $b_0 \neq 0$. Let $s'=b_1r+b_0s$ and $r'=r$. By substituting $s=(s'-b_1r')/b_0$ and $r=r'$ into $f$ and, if necessary, multiplying $f$ by a square to make it integral again we can assume that the linear factor is $s$. Let $f'$ denote the result of this substitution and let $a_0',\ldots, a_4'$ denote it's coefficients.  Then $a_4'=0$ because $s|f'$, and $a_3'\neq 0$ since we assumed that $f$ and hence $f'$ doesn't have any repreated roots.
	
	Consider the resulting equation $$y'^2= s'(a'_3r'^3+a'_2r'^2s'+a'_1r's'^2+a'_0s'^3)$$  then any prime dividing $s$ with an odd valuation also needs to divide $(a'_3r'^3+a'_2r'^2s'+a'_1r's'^2+a'_0s'^3)$ and hence $a'_3r'^3=a'_3(b_1r+b_0s)^3$. The coprimality of $r$ and $s$ ensures that such a prime needs to divide $a'_3b_1$. Let $S$ be a finite set of primes of $\OO_K$ that generate the class group of $K$ that also contains the divisors of $a'_3b_1$. Then $s' \in K^*/K^{*2}$ lands in the finite subgroup $H \subseteq  K^*/K^{*2}$ generated by  $\OO_{K,S}^*$, the $S$-units in $\OO_K$. In particular for every solution $x,y,z$ of $y^2=f(x^2,z)$ with $x,z$ coprime there is a $d \in H$ such that $x,y,z$ comes from a point on the hyperelliptic genus 2 curve $y^2=f'(x^2,dw^2)$. And genus 2 curves only have finitely many rational points.
\end{proof}  

\begin{remark}There are algorithms to compute the ring of $S$-units $\OO_{K,S}^*$. So the proof of the above theorem is also gives an strategy for finding all solutions to the above equation, given that one can find the rational points on the genus 2 curves involved.
\end{remark}
	
	\section{Divisibility by 2}	\label{sec:divisibility_by_2}

	Genus field theory as developed by Gauss \cite{gauss66disquisitiones} describes the $2$ torsion in the narrow class group of quadratic fields. Since for imaginary quadratic fields the class group and the narrow class group are the same, this allows us to easily determine the divisibility by $2$ of class groups of imaginary quadratic fields. This is made precise in the following theorem 
	\begin{theorem}\label{thm:genus_field}
		Let $n\geq 1$ be and integer,  $p_1,\ldots, p_n$ be distinct primes and $d = \pm \prod_{i=1}^n p_i$ be a squarefree integer. Then the $2$-rank of the narrow class group of $\Q(\sqrt{d})$ is $n-1$. If $d<0$ then the 2-rank of the class group of $\Q(\sqrt{d})$ is also $n-1$.
	\end{theorem}

With this result in hand we can now start the proof of the divisibility by $2$ part.

	\begin{theorem}\label{thm:divisibility_by_2}
	Let $K$ be an imaginary quadratic field such that $Y_1(16)(K) \neq \emptyset$   then the class number of $K$ is even. 
\end{theorem}

\begin{proof}

In the introduction an equation for $X_1(16)$ is already given in \cref{eq:x1_16}. However for this proof it is more convenient to work with the homogenized version of this equation 	\begin{align}
\begin{split}
X_1(16) : t^2 &= h_{16}(r,s) \\
h_{16}(r,s) &:= rs(r^2 + s^2)(r^2 + 2 rs - s^2) \label{eq:x1_16_h}
\end{split}
\end{align}
in the weighted projective space $\PP(1,3,1)$. Here $r,s$ have weight 1 and $t$ has weight 3. The map to the affine model in \cref{eq:x1_16} is given by $(r:t:s) \mapsto (r/s,t/s^3)$.

If $r,s$ are coprime integers with $s\neq 0$ then points with $x=r/s$ are defined over the field $\Q(\sqrt{h_{16}(r,s)})$. It is already known that $Y_1(16)(\Q(-1)) = \emptyset$ \cite[Thm. 1]{najman2010cyclotomicquadratictorsion}. Furthermore, the only solutions with $t=0$ are cusps, so this means that \Cref{thm:genus_field} reduces the proof of \Cref{thm:divisibility_by_2} to part $(a)$ of the following proposition.

\end{proof}

\begin{proposition}\label{prop:divisibilit_by_2}\,
	
\begin{itemize}
\item[a)] The equation $-p t^2 = h_{16}(r,s)$ has no solutions with $p,r,s,t\in \mathbb Z$ , $p>0$ a prime $r,s$ coprime and $t\neq 0$. 
\item[b)] Let $p,r,s,t\in \mathbb Z$ , with $p>0$ a prime $r,s$ coprime and $t\neq 0$ be a solution  to $p t^2 = h_{16}(r,s)$. After a transformation of the form $(r,s) \mapsto (\pm r, \pm s)$ or  $(r,s) \mapsto (\pm s, \mp r)$ there are integers $u,v$ and $y,z$ such that $r=u^2$, $s=v^2$, $2 h_2(u^2,v^2)=y^2$ and $2 h_1(u^2,v^2)=pz^2$.
\end{itemize}
\end{proposition}

The proof below is done by a repeated $2$ descent. This two descent allows us to relate the solutions to  $\pm p t^2 = h_{16}(r,s)$ to rational points on curves a finite set of auxiliary curves that do not depend on the prime $p$. The proposition is then proved by determining the rational points on all curves in this finite set. The proof is rather lengthy since the two descent leads to a lot of different case distinctions. The cases and the running assumptions have been clearly marked in bold. Cases are split up into further sub-cases on an as needed basis. For example 2.a.i is a sub-case of 2.a which is a sub-case of case 2. 
\begin{proof}
By multiplying $r$ and $s$ by $-1$ if necessary we can reduce to the case $r\geq 0$. 	
	
Note that $r,s$ and $(r^2+s^2)(r^2+2 rs - s^2)$ are pairwise coprime. So for $h_{16}(r,s) := rs(r^2+s^2)(r^2+2 rs - s^2)$ to equal $\pm pt^2$  one needs at least $2$ of the three values $r,s$ and $(r^2+s^2)(r^2+2 rs - s^2)$ to be $\pm$ a square. Note that $h_{16}(r,s) = h_{16}(s,-r)$, so we can reduce to the case that $r$ is a square and that $u$ is an integer such that $r=u^2$.

{\bf From now on $\mathbf{r=u^2}$.}

Let $h_1(r,s) := (r^2+s^2)$ and $h_2(r,s) := (r^2+2 rs - s^2)$. Then one has
\begin{align*} 
4r^3 &= (3r-s)h_1+(r-s)h_2 \\
4s^3 &=  (3s+r)h_1 -(r+s)h_2. 
\end{align*}
So the coprimality of $r,s$ implies $\gcd(h_1(r,s),h_2(r,s)) \mid 4$.
The proof is broken up into two cases, depending on whether $\pm s$ or $\pm (r^2+s^2)(r^2+2 rs - s^2)=\pm h_1(r,s)h_2(r,s)$  is a square.

{\bf Case 1; $\mathbf{s = \pm v^2}$ :} 
The fact that $\gcd(h_1(r,s),h_2(r,s)) \mid 4$ implies there is an $k\in 0,1$ and $y\in \Z$ such that either $\pm 2^ih_1(r,s)=y^2$ or $\pm 2^ih_2(r,s)=iy^2$. This gives rise to another case distinction.

{\bf Case 1.a; $\mathbf{\pm 2^kh_1(r,s)=y^2}$:}
Note that $h_1(r,s) \geq 0$ so that $2^kh_1(r,s)=y^2$. This means we are reduced to finding the rational points on the two elliptic curves in $\PP(1,2,1)$:

\begin{align*}
E_1 : y^2 &= (u^4+v^4) \\ 
E_2 : y^2 &=2(u^4+v^4)
\end{align*}
The elliptic curves both have rank 0 \github{order_2.m}{Claim 1} and torsion subgroups isomorphic to $(\Z/2\Z)^2$ \github{order_2.m}{Claim 2}. The points on $E_1$ satisfy $(r:s)\in \set{(0:1), (1:0)}$ and the points on $E_2$ satisfy  $(r:s)\in \set{(-1:1), (1:1)}$  \github{order_2.m}{Claim 3}. Note that $(0:1)$ and $(1:0)$ both give $t=0$. Since $h_{16}(-1,1) = h_{16}(1,1)=2^2$ the points $(-1:1)$ and $(1:1)$ give $t=2$ and $p=1$. So none of the points lead to a valid solution of the original equation. This also proves that $p t^2 = h_{16}(r,s)$ has no solutions in case $1a$.

{\bf Case 1.b; $\mathbf{\pm 2^kh_2(r,s)=y^2}$:}
Note that $h_2(r,s)=-h_2(s,-r)=h_2(-s,r)$, so we can restrict to  $2^kh_2(r,s)=y^2$ if we later remember to also add $(s: -r)=(-s: r)$ for every point $(r:s)$ that we have found. This means we need to find all rational points on the following 4 genus 1 curves in $\PP(1,2,1)$.

\begin{align}
\begin{split}
E_3 : y^2 &= (u^4+2u^2v^2-v^4) \\ 
E_4 : y^2 &=2(u^4+2u^2v^2-v^4)\label{eq:E4} \\
E_5 : y^2 &= (u^4-2u^2v^2-v^4) \\ 
E_6 : y^2 &=2(u^4-2u^2v^2-v^4).
\end{split}
\end{align}

The curve $E_6$ doesn't have any rational points since the lefthand side is  $0, 1, 4$ or $9  \mod 16$ while the right hand side is  $2, 12$ or $14  \mod 16$ whenever $2 \nmid gcd(u,v)$ \github{order_2.m}{Claim 4}.

For the 3 remaining curves note that $(1:1:0)$,$(1:2:1)$ and $(1:1:0)$ are points on $E_3, E_4$ and $E_5$ respectively, turning them into elliptic curves. These curves have rank $0, 1$ and $0$ respectively \github{order_2.m}{Claim 5} and they all have torsion subgroup isomorphic to $\Z/2\Z$ \github{order_2.m}{Claim 6}.

For $E_3$ and $E_5$ one can easily compute that points on them have to lead to $(r:s)=(1:0)$ \github{order_2.m}{Claim 7}. This point and $(s:-r) = (0:1)$ both give $t=0$ so are not valid solutions to the original equation.

For points on $E_4$ we note that this case corresponds to  $s=v^2$. 
This means that $$h_{16}(r,s)=u^2v^2(u^4+v^4)(u^4+2u^2v^2-v^4))=\frac 1 2 u^2v^2(u^4+v^4)y^2 >0.$$
In particular $h_{16}(r,s) \neq -pt^2$.
The argument here also show that in case 1 the only solutions to $pt^2=h_{16}(r,s)$ are coming from points $(u,y,v)$ on $E_4$ such that $2(u^4+u^4)=2h_1(r,s)$ is of the form $pz^2$ for some integer $z$. 

{\bf Case 2; $\mathbf{\pm h_1(r,s) h_2(r,s)=y^2}$:}
If the following two genus 1 curves were elliptic curves of rank 0, then this would be easily solvable.
\begin{align}
\begin{split}
E_7 : y^2 &= h_1(r,s) h_2(r,s) \label{eq:E7}\\ 
E_8 : y^2 &= -h_1(r,s) h_2(r,s) \\
\end{split}
\end{align}
However, both these elliptic curves have rank 1 \github{order_2.m}{Claim 8}, so more work is needed. However, since $r=u^2$ we already know there will only be finitely many solutions due to \Cref{thm:diophantine_mainv2}. We will mimic the finiteness proof of \Cref{thm:diophantine_mainv2} to determine all solutions in this case, with one modification. Instead of going to a quadratic field extension to find a linear factor of $ h_1(r,s) h_2(r,s)$, the proof will use the factorisation into two quadratic factors over $\Z$.

Since $\gcd(h_1(r,s), h_2(r,s)) \mid 4$ and $h_1(r,s)\geq 0$. This means that there is an integer $z$ such that either $h_1(r,s) = z^2$ or $h_1(r,s) = 2z^2$, giving another case distinction.

{\bf Case 2.a;  $\mathbf{h_1(r,s) = z^2}$:}
In this case $r^2+s^2=z^2$. So Euclid's formula for Pythagorean triples ensures that there are coprime integers $n,m$ such that either 
\begin{itemize}
\item[i)] $r = m^2-n^2, s=2mn, z=m^2+n^2$ or,
\item[ii)] $r = 2mn, s = m^2-n^2, z=m^2+n^2$.
\end{itemize}
These two cases are dealt with separately.

{\bf Case 2.a.i; $\mathbf{r = m^2-n^2, s=2mn}$:}
In this case we have $u^2+n^2=m^2$. So there are rational numbers $k$ and $v$ such that $$u=k(v^2-1), n=2kv, m=k(v^2+1).$$ Substituting everything back gives 
$$y^2=\pm z^2h_2(r,s)=\pm z^2k^4(v^8 + 8v^7 - 20v^6 - 8v^5 - 26v^4 - 8v^3 - 20v^2 + 8v + 1)$$ 
as  verified in \github{order_2.m}{Claim 9}.
 In particular any solution in this case will give a rational point on one of the two genus 3 curves below:

\begin{align*}
\begin{split}
C_9: y'^2 &= (v^8 + 8v^7 - 20v^6 - 8v^5 - 26v^4 - 8v^3 - 20v^2 + 8v + 1) \label{eq:C1}\\ 
C_{10} : y'^2 &= -(v^8 + 8v^7 - 20v^6 - 8v^5 - 26v^4 - 8v^3 - 20v^2 + 8v + 1)\\
\end{split}
\end{align*}

Both these curves admit a map of degree 2 to 
\begin{align*}
D_{10}: y^2 = x^5 + 4x^4 - 6x^3 - 4x^2 + x  \quad\quad \text{\github{order_2.m}{Claim 10}}.
\end{align*}

The Mordell-Weil group of $J(D_{10})$ is isomorphic to $\Z/2\Z \times \Z$  \github{order_2.m}{Claim 11}, so it is of rank 1.  Applying Chabauty's method to determine all the rational points on $D_{10}$ shows that it has 4 rational points, namely $(0,0)$, $(-1,-2)$, $(-1,2)$ and the point at infinity \github{order_2.m}{Claim 12}. The rational points on $C_9$ and $C_{10}$ can be computing by looking at the inverse images under $C_9 \to D_{10}$ and $C_{10} \to D_{10}$ of these 4 points. In doing so one finds that $(0,-1)$ and $(0,1)$ together with two points at infinity are the only rational points on $C_9$ and that the 4 points of the form $(\pm 1, \pm 8)$ are the only rational points on $C_{10}$.

Transferring the points on $C_{9}$ and $C_{10}$ back to the original equation show that the rational points on $C_9$ all lead to $(r:s)=(1:0)$ while those on $C_{10}$ lead to $(r:s)=(0:1)$. As mentioned before, both of these are not solutions to the original equation.

{\bf Case 2.a.ii: $\mathbf{r = 2mn, s = m^2-n^2}$:} In this case $u^2=2nm$. In particular, there exists integers $v,w$ such that either 
\begin{enumerate}
	\item[A)]$n=\pm v^2$ and $m=\pm 2w^2$ or 
	\item[B)]$n=\pm 2v^2$ and $m=\pm w^2$. 
\end{enumerate}	
Since $m$ and $n$ have the same sign, and changing the signs of both $m$ and $n$ at the same time doesn't change $r$ and $s$ we can assume that both are positive.

Setting $y'=y/z$ substituting everything back and de-homogenizing by setting $w=1$ gives the curves
\begin{align*}
	\begin{split}
		C_{11}: y'^2 &= v^8 - 8v^6 - 24v^4 + 32v^2 + 16 \\ 
		C_{12} : y'^2 &= -(v^8 - 8v^6 - 24v^4 + 32v^2 + 16)\\
	\end{split}
\end{align*}
 in case A) and  the curves
 \begin{align*}
 	\begin{split}
 		C_{13}: y'^2 &= 16v^8 - 32v^6 - 24v^4 + 8v^2 + 1\\ 
 		C_{14} : y'^2 &= -(16v^8 - 32v^6 - 24v^4 + 8v^2 + 1)\\
 	\end{split}
 \end{align*}
 in case B) \github{order_2.m}{Claim 15} . The map $(y',v')\mapsto (1/4y',1/2v)$ induces an isomorphism between  $C_{11}$ and $C_{13}$ as well as an isomorphism betweeen $C_{12}$ and $C_{14}$. So it suffices to determine the rational points on $C_{11}$ and $C_{12}$.

The curves $C_{11}$ and $C_{12}$ both admit a map of degree 2 to the curve
\begin{align*}
	D_{12}: y^2 =  2x^5 - 8x^4 - 12x^3 + 8x^2 + 2x  \quad\quad \text{\github{order_2.m}{Claim 16}}.
\end{align*}

The Mordell-Weil group of $J(D_{12})$ has rank 0 \github{order_2.m}{Claim 17}.  Applying Chabauty's method to determine all the rational points on $D_{12}$ shows that it has 2 rational points, namely $(0,0)$ and the point at infinity \github{order_2.m}{Claim 18}. Checking the points on $C_{11}$ and $C_{12}$ mapping to these two points shows that $C_{11}$ has exactly four rational points, namely $(0,\pm 4)$ and the two points at infinity \github{order_2.m}{Claim 19} .  These four points on $C_{11}$ all lead to $r=0$, so don't correspond to solutions of the original problem.  The two rational points on $D_{12}$ both pull back to degree 2 points on $C_{12}$ \github{order_2.m}{Claim 20}, implying that $C_{12}$ and hence $C_{14}$ don't have any rational points. The 4 rational points on $C_{13}$ also translate to solutions with $r=0$, so don't correspond to solutions of the original problem.

{\bf Case 2.b; $\mathbf{h_1(r,s) = 2z^2}$:}
Note that $(1:1:1)$ is a point on the quadric $h_1(r,s) = r^2+s^2 = 2z^2$. By considering lines trough this point, one can show that up to scaling every point on this quadric is of the form $$(m^2 - 2mn - n^2)^2 + (m^2 +2mn - n^2)^2 = 2(m^2+n^2)^2$$ for coprime integers $m,n$.


Let $f_1 := m^2 - 2mn - n^2$ and $f_2 := m^2 + 2mn - n^2$, then $(2m-n)f_1 + (2m+n)f_2=4m^3$ and $(-m-2n)f_1  + (m-2n)f_2=4n^3$ so that 
$$\gcd(f_1, f_2) \mid \gcd(4m^3,4n^3) =4.$$
By considering $f_1$ mod $4$ one can show that $4 \nmid f_1$ for coprime $m,n$ so that $\gcd(f_1, f_2) \mid 2$.


In particular there are two cases to consider
\begin{enumerate}
	\item[i)] $r = m^2 - 2mn - n^2$ and $s = m^2 + 2mn - n^2$,
	\item[ii)] $r = (m^2 - 2mn - n^2)/2$ and $s = (m^2 + 2mn - n^2)/2$. 
\end{enumerate}	
Note in case ii one needs $m$ and $n$ to both be odd.

{\bf Case 2.b.i; $\mathbf{r = m^2 - 2mn - n^2}$ and $\mathbf{s = m^2 + 2mn - n^2}$:}

In this case $u^2 = m^2 - 2mn - n^2$, which can be rewritten as $u^2+(m+n)^2=2m^2$.  This means there exists integers $v,w$ and a rational number $k$ such that $u =k(v^2 - 2vw - w^2)$,  $n = k(2vw - 2w^2)$.

Substituting all of this into $y^2=\pm h_1(r,s)h_2(r,s)$ and dehomogenizing by setting $w=1$ we get that we need to find the rational solutions to
$$y^2=\pm 4z^2k^4(v^8 - 8v^7 - 12v^6 + 56v^5 - 122v^4 + 136v^3 - 76v^2 + 72v - 31)$$
as verified in \github{order_2.m}{Claim 21}. The factor $4z^2k^4=(2zk^2)^2$ can be absorbed into $y$, leading to the curves
\begin{align*}
	\begin{split}
		C_{15}: y'^2 &= v^8 - 8v^7 - 12v^6 + 56v^5 - 122v^4 + 136v^3 - 76v^2 + 72v - 31\\ 
		C_{16} : y'^2 &= -(v^8 - 8v^7 - 12v^6 + 56v^5 - 122v^4 + 136v^3 - 76v^2 + 72v - 31)\\
	\end{split}
\end{align*}

Both of these curves have a degree 2 map to:
\begin{align*}
	 D_{16}: y^2 &= x^6 - 6x^5 + x^4 + 20x^3 - x^2 - 6x - 1.  \quad\quad \text{\github{order_2.m}{Claim 22}}
\end{align*}

The Mordell-Weil group of $J(D_{16})$ is isomorphic to $\Z/2\Z \times \Z$  \github{order_2.m}{Claim 23}, so it is of rank 1.  Applying Chabauty's method to determine all the rational points on $D_{16}$ shows that it only has the two points at infinity \github{order_2.m}{Claim 24} as rational points. The rational points on $C_{15}$ and $C_{16}$ can be computing by looking at the inverse images of these 4 points. In doing so one finds that $(1,\pm 4)$ together with two points at infinity are the only rational points on $C_{15}$ \github{order_2.m}{Claim 25} and that $C_{16}$ has no rational points \github{order_2.m}{Claim 26}.

Transferring the points on $C_{15}$ back to the original equation show that the rational points on $C_{15}$
all lead to solutions with $(r:s)=(1:1)$ and $p=1$, which are not solutions to the original equation.

{\bf Case 2.b.ii; $\mathbf{r = (m^2 - 2mn - n^2)/2}$ and $\mathbf{s = (m^2 + 2mn - n^2)/2}$:}

In this case $u^2 = (m^2 - 2mn - n^2)/2$, which can be rewritten as $u^2+n^2=2((m-n)/2)^2$.  This means there exists integers $v,w$ and a rational number $k$ such that $u =k(v^2 - 2vw - w^2)$,  $n = k(v^2+2vw - 2w^2)$.

Substituting all of this into $y^2=\pm h_1(r,s)h_2(r,s)$ and dehomogenizing by setting $w=1$ we get that we need to find the rational solutions to
$$y^2=\pm 4z^2k^4(-17v^8 - 104v^7 - 132v^6 - 72v^5 - 22v^4 + 40v^3 + 28v^2 + 8v - 1)$$
as verified in \github{order_2.m}{Claim 27}. The factor $4z^2k^4=(2zk^2)^2$ can be absorbed into $y$, leading to the curves
\begin{align*}
	\begin{split}
		C_{17}: y'^2 &= -17v^8 - 104v^7 - 132v^6 - 72v^5 - 22v^4 + 40v^3 + 28v^2 + 8v - 1\\ 
		C_{18} : y'^2 &= -(-17v^8 - 104v^7 - 132v^6 - 72v^5 - 22v^4 + 40v^3 + 28v^2 + 8v - 1)\\
	\end{split}
\end{align*}

Both of these curves have a degree 2 map to:
\begin{align*}
	D_{18}: y^2 &= x^6 - 2x^5 - 15x^4 - 4x^3 + 15x^2 - 2x- 1.  \quad\quad \text{\github{order_2.m}{Claim 28}}
\end{align*}

The Mordell-Weil group of $J(D_{18})$ is isomorphic to $\Z/2\Z \times \Z$  \github{order_2.m}{Claim 29}, so it is of rank 1.  Applying Chabauty's method to determine all the rational points on $D_{18}$ shows that it only has the two points at infinity \github{order_2.m}{Claim 30} as rational points. The rational points on $C_{17}$ and $C_{18}$ can be computing by looking at the inverse images of these 4 points. In doing so one finds that  $C_{17}$ has no rational points \github{order_2.m}{Claim 31} and that $(0,\pm 1)$ and $(-1,\pm 4)$ are the only rational points on $C_{18}$ \github{order_2.m}{Claim 32}.

Transferring the points on $C_{18}$ back to the original equation show that the rational points on $C_{18}$
all lead to solutions with $(r:s)=(1:-1)$ and $p=1$, which are not solutions to the original equation.

\end{proof}

\section{Narrow class number of real quadratic points}\label{sec:real_quadratic_points}

\begin{question}\label{Q:realquadratic1}
	Is the set of real quadratic fields $K$ with odd narrow class number such that $Y_1(16)(K) \neq \emptyset$ finite? Or equivalently, are there only finitely many $p,r,s,t\in \mathbb Z$ with $p$ a prime, $r,s$ coprime $t \neq 0$ such that $p t^2 = h_{16}(r,s)$. 
\end{question}

By \Cref{prop:divisibilit_by_2} the above question is equivalent to the following:
\begin{question}\label{Q:realquadratic2}
	Are there infinitely many rational points on the rank 1 elliptic curve $E_4$ from \cref{eq:E4} such that $2(u^4+v^4)$ is of the form $pz^2$ for some integer $z$ and a prime $p$?
\end{question}

\begin{remark}The quantity $2(u^4+v^4)$ depends on the projective coordinates chosen to represent the rational point on $E_4$. However choosing different coordinates cause it to change by a 4-th power, so that the question whether it is of the form $pz^2$ is a well defined question independent of the representative.
\end{remark}

\begin{example} By looking for points on $E_4$ such that $2(u^4+v^4)$ is of the form $pz^2$, one can find some solutions to $p t^2 = h_{16}(r,s)$. The largest one found was the following. 
	Define the integers:
	\begin{align*}
	u :=& 1383308224231610113228232741369733180270315041 \\
	v:=& 702229330665242264680897734882798122886724801\\
	y :=&2\cdot 162875761464289218878723159174898919197655635215733128961834\\ &3575391414712185591538941481919 \\
	z :=& 2\\
	p :=& h_1(u^2,v^2)/2\\
	r :=& u^2 \\
	s :=& v^2 
	\end{align*}
	then $p$ is a prime of $181$ digits and
	\begin{align*}
	(y)^2 &=2(u^4+2u^2v^2-v^4)\\
	 p(yuv)^2&= h_{16}(r,s).
	\end{align*}
	So contrary to the imaginary quadratic case, there is a large example of a point in $Y_1(16)(K)$ where $K$ is a real quadratic field with odd narrow class number.
\end{example}

\subsection{Heuristic answer to Questions \ref{Q:realquadratic1} and \ref{Q:realquadratic2}. }
In this section we will discuss a heuristic that suggests that the answer to Questions $6.1$ and $6.2$ will be no. Implying that there will be only finitely quadratic fields $K$ with odd narrow class number such that $Y_1(16)(K)\neq \emptyset$.

Let $n$ be an integers and denote by $\pi_2(n)$ the number of positive integers smaller then $n$ of the form $p z^2$ with $p$ prime and $z$ and integer. By considering the case $z=1$ it is clear that $\pi_2(n)$ grows as least as fast as $c_1 n/\log(n)$ for some constant $c_1>0$.  The following 
rough estimate shows that similary $\pi_2(n)$ is bounded from above by $c_2 n/\log(n)$ for some different constant. 

\begin{align*}
\pi_2(n) \sim& \sum_{k=1}^{\lfloor\sqrt{n/2}\rfloor} \frac n {k^2\log(n/k^2)} \\
=&
 n\left(\sum_{k=1}^{\lfloor\sqrt[4]{n}\rfloor} \frac 1 {k^2\log(n/k^2)}  +
\sum_{k={\lfloor\sqrt[4]{n}\rfloor}+1}^{\lfloor\sqrt{n/2}\rfloor} \frac 1 {k^2\log(n/k^2)}  \right) \\
\leq &  n\left(\sum_{k=1}^{\lfloor\sqrt[4]{n}\rfloor} \frac 2 {k^2\log(n)}  +
\sum_{k={\lfloor\sqrt[4]{n}\rfloor}+1}^{\lfloor\sqrt{n/2}\rfloor} \frac 1 {k^2\log(2)}  \right) \\
\leq &  n\left(\frac {\pi^2} {3 \log(n)}+
O(\frac{1}{\sqrt[4]{n}})\right) =O\left(\frac{n} {\log(n)}\right).
\end{align*}

This suggests a heuristic model that predicts that an integer of size about $n$ is of the form $pz^2$ with a chance of $c/\log(n)$ for some constant $n$. Now $E_4$ is of rank $1$ and hence (up to torsion) generated by a single point $P$. Now if $m$ is an integer then the logarithmic height of $mP$ grows as $m^2h(P)$. In particular $log(2(u(mP)^4+v(mP)^4))$ will grow as $c'm^2$ for some constant $c'$. As a conceauence we expect there to be roughly $\sum_{m=1}^\infty 2\#E_4(\Q)_{tors} c/c' m^{-2}$ points in $E_4(\Q)$ for which $2(u^4+v^4)$ is of the form $pz^2$. Since that sum converges the heuristic suggest that there should only be finitely many points in $E_4(\Q)$ for which $2(u^4+v^4)$ is of the form $pz^2$.
	
	\printbibliography

@article{bruin2009two,
  title={Two-cover descent on hyperelliptic curves},
  author={Bruin, Nils and Stoll, Michael},
  journal={Mathematics of computation},
  volume={78},
  number={268},
  pages={2347--2370},
  year={2009}
}

@article{bruin2003chabauty,
url = {https://doi.org/10.1515/crll.2003.076},
title = {Chabauty methods using elliptic curves},
author = {Nils Bruin},
pages = {27--49},
volume = {2003},
number = {562},
journal = {Journal für die reine und angewandte Mathematik},
doi = {doi:10.1515/crll.2003.076},
year = {2003},
}

@article {magma,
    AUTHOR = {Bosma, Wieb and Cannon, John and Playoust, Catherine},
     TITLE = {The {M}agma algebra system. {I}. {T}he user language},
      NOTE = {Computational algebra and number theory (London, 1993)},
   JOURNAL = {J. Symbolic Comput.},
  FJOURNAL = {Journal of Symbolic Computation},
    VOLUME = {24},
      YEAR = {1997},
    NUMBER = {3-4},
     PAGES = {235--265},
      ISSN = {0747-7171},
   MRCLASS = {68Q40},
       DOI = {10.1006/jsco.1996.0125},
       URL = {http://dx.doi.org/10.1006/jsco.1996.0125},
}

@phdthesis{krummthesis,
  TITLE = {Quadratic points on modular curves},
  AUTHOR = {Krumm, David},
  SCHOOL = {University of Georgia},
  YEAR = {2013},
  TYPE = {Thesis},
}

@article{najman2010cyclotomicquadratictorsion,
	author = {Najman, Filip},
	title = {Complete classification of torsion of elliptic curves over quadratic cyclotomic fields},
	fjournal = {Journal of Number Theory},
	journal = {J. Number Theory},
	issn = {0022-314X},
	volume = {130},
	number = {9},
	pages = {1964--1968},
	year = {2010},
	language = {English},
	doi = {10.1016/j.jnt.2009.12.008},
	zbMATH = {5759066},
	Zbl = {1200.11039}
}

@article{mazur_mod_curves,
  title={Modular curves and the {E}isenstein ideal},
  author={Mazur, Barry},
  journal={Publications Math{\'e}matiques de l'Institut des Hautes {\'E}tudes Scientifiques},
  volume={47},
  number={1},
  pages={33--186},
  year={1977},
  publisher={Springer},
  note={With an appendix by Barry Mazur and Michael Rapoport}
}

@Book{liu_geom_curves,
	Author = {Liu, Qing},
	Title = {Algebraic geometry and arithmetic curves. {Transl}. by {Reinie} {Ern{\'e}}},
	FSeries = {Oxford Graduate Texts in Mathematics},
	Series = {Oxf. Grad. Texts Math.},
	Volume = {6},
	ISBN = {0-19-920249-4},
	Year = {2006},
	Publisher = {Oxford: Oxford University Press},
	Language = {English},
	zbMATH = {5048200},
	Zbl = {1103.14001}
}

@misc{banwait_derickx_quadratic_torsion,
	title={Torsion subgroups of elliptic curves over quadratic fields and a conjecture of Granville}, 
	author={Barinder S. Banwait and Maarten Derickx},
	year={2024},
	eprint={2401.14514},
	archivePrefix={arXiv},
	primaryClass={math.NT},
	url={https://arxiv.org/abs/2401.14514}, 
}

@misc{gauss66disquisitiones,
	Author = {Gauss, Carl Friedrich},
	Title = {Disquisitiones arithmeticae. ({Translated} by {Arthur} {A}. {Clarke})},
	Year = {1966},
	Language = {English},
	HowPublished = {New {Haven}-{London}: {Yale} {University} {Press}. xx, 472 p. (1966).},
	zbMATH = {3221502},
	Zbl = {0136.32301}
}

@article{GillibertLevenTorsiontoIdealClasses,
	author = {Gillibert, Jean and Levin, Aaron},
	title = {Pulling back torsion line bundles to ideal classes},
	fjournal = {Mathematical Research Letters},
	journal = {Math. Res. Lett.},
	issn = {1073-2780},
	volume = {19},
	number = {6},
	pages = {1171--1184},
	year = {2012},
	language = {English},
	doi = {10.4310/MRL.2012.v19.n6.a1},
	zbMATH = {6249184},
	Zbl = {1325.14019}
}
	
\end{document}